\documentclass[12pt,a4paper]{article}
\UseRawInputEncoding 
\usepackage{amstext}
\usepackage{fancyhdr}
\usepackage{amsfonts,graphicx,bezier, amssymb}
\usepackage{amsmath}
\usepackage{caption}
\usepackage{mathtools}
\usepackage{tikz}
\usepackage{rotating}

\usetikzlibrary{arrows.meta}
\newenvironment{prf}{\noindent{\bf{Proof:}}~~}{\hfill\rule{1ex}{1ex}\vskip1.5ex}
\newcommand{\Z}{\mathbb Z}
\newcommand{\Q}{\mathbb Q}

\newcommand{\beqa}{\begin{eqnarray}}
\newcommand{\enqa}{\end{eqnarray}}
\newcommand{\beq}{\begin{eqnarray*}}
\newcommand{\enq}{\end{eqnarray*}}

\newtheorem{rem}{Remark}[section]

\newtheorem{cor}{Corollary}[section]
\newtheorem{propn}{Proposition}[section]
\newtheorem{defn}{Definition}[section]
\newtheorem{exam}{{\bf Example}}[section]
\newtheorem{thm}{Theorem}[section]

\newtheorem{lem}{Lemma}[section]

\newcommand{\noi}{\noindent}

\makeatletter
\providecommand*{\twoheadrightarrowfill@}{%
  \arrowfill@\relbar\relbar\twoheadrightarrow
}
\providecommand*{\twoheadleftarrowfill@}{%
  \arrowfill@\twoheadleftarrow\relbar\relbar
}
\providecommand*{\xtwoheadrightarrow}[2][]{%
  \ext@arrow 0579\twoheadrightarrowfill@{#1}{#2}%
}
\providecommand*{\xtwoheadleftarrow}[2][]{%
  \ext@arrow 5097\twoheadleftarrowfill@{#1}{#2}%
}
\makeatother

\begin{document}

\begin{center}
{\bf\Large Applications of reduced and coreduced modules I}
\end{center}

\vspace*{0.3cm}
\begin{center}
David Ssevviiri\\
\vspace*{0.3cm}
Department of Mathematics\\
Makerere University, P.O BOX 7062, Kampala Uganda\\
Email: david.ssevviiri@mak.ac.ug, ssevviiridaudi@gmail.com
\end{center}

\begin{abstract}This is the first in a series of papers  highlighting  the applications of reduced and coreduced modules. Let $R$ be a commutative unital ring and $I$ an ideal of $R$. 
We show that   $I$-reduced $R$-modules and   $I$-coreduced $R$-modules provide a setting  in which the Matlis-Greenless-May (MGM) Equivalence and the Greenless-May (GM) Duality hold. These two notions have been hitherto only known to exist in the derived category setting. We realise the $I$-torsion and the $I$-adic completion functors as representable functors and under suitable conditions compute natural transformations between them and other functors.
\end{abstract}


\vspace*{0.4cm}

{\bf MSC 2010} Mathematics Subject Classification: 13D07,  13C12, 13B35, 13C13 

 \section{Introduction}

 \begin{paragraph}\noi 
Let $R$ be a commutative unital ring and $I$ an ideal of $R$.
The additive endo functors $
\Gamma_I ~:~R\text{-Mod}\rightarrow R\text{-Mod}$, $M\mapsto \Gamma_I(M):= \lim\limits_{\rightarrow}\text{Hom}_R(R/I^k, M)$ and $ \Lambda_I ~:~R\text{-Mod}\rightarrow R\text{-Mod}$,
$M\mapsto \Lambda_I(M):=\lim\limits_{\leftarrow} (M/I^kM)$; called the {\it $I$-torsion functor} and the {\it $I$-adic completion functor} respectively have been widely studied. Grothendieck was the first to study the torsion functor in the algebraic geometry setting of sheaves, \cite{Hartshorne}. This functor is central to the well known notions of local cohomology, local duality, \cite{Brodmann}; Cousin complexes, \cite{Kempf}; and their generalisations, \cite{glduality, Suzuki, Yassemi}. Sharp in \cite{Sharp} introduced their algebraic avatars. $\Lambda_I$ the dual to the functor $\Gamma_I$   has been used widely to study local homology, see \cite{Divaani,  Greenless, Peter} among others. The functors $\Gamma_I$ and  $\Lambda_I$ and their derived functors were key in \cite{Alonso, Barthel, Barthel2, Greenless, Yekutieli, Peter, Y2, Vyas} where notions such as Greenless-May (GM) Duality and Matlis-Greenless-May (MGM) Equivalence were studied in different settings. \end{paragraph}

\begin{paragraph}\noi 
 The functors $\Gamma_I$ and $\Lambda_I$ are not adjoint to each other. However, under suitable conditions, their corresponding derived functors $${\bf R}\Gamma_I~~\text{and}~~ {\bf L}\Lambda_I: D(R)\rightarrow D(R)$$ where $D(R)$ is the derived category of $R$-modules are adjoint. This is what is known as the Greenless-May (GM) Duality.   Our first main result, is Theorem \ref{01} which shows that in the setting of $I$-reduced modules and $I$-coreduced modules, the functors $\Gamma_I$ and $\Lambda_I$ form an adjoint pair.
\end{paragraph}

\begin{paragraph}\noi 
 Let $I$ be an ideal of a ring $R$ and   $(R\text{-Mod})_{I\text{-red}}$  (resp. $(R\text{-Mod})_{I\text{-cor}}$) denote a full subcategory of $R$-Mod consisting of all $I$-reduced (resp. $I$-coreduced) $R$-modules. 
 \end{paragraph}

\begin{thm}\label{01}{\rm {\bf [GM Duality in $R$-Mod]}}
 For any ideal $I$ of a ring $R$, 
 \begin{enumerate}
  \item  The functor 
 $$\Gamma_I: (R\text{-Mod})_{I\text{-red}} \rightarrow (R\text{-Mod})_{I\text{-cor}}$$ is idempotent and for any $M\in (R\text{-Mod})_{I\text{-red}}$, $\Gamma_I(M)\cong \text{Hom}_R(R/I, M)$.
 
 \item   The functor $$\Lambda_I: (R\text{-Mod})_{I\text{-cor}} \rightarrow (R\text{-Mod})_{I\text{-red}}$$ is idempotent and  for any  $M\in (R\text{-Mod})_{I\text{-cor}}$, $\Lambda_I(M)\cong R/I\otimes M$.
 \item  For any $N\in (R\text{-Mod})_{I\text{-red}}$ and $M\in (R\text{-Mod})_{I\text{-cor}}$, 
 
 $$\text{Hom}_R(\Lambda_I(M), N)\cong \text{Hom}_R(M, \Gamma_I(N)).$$
 \end{enumerate}

\end{thm}

\begin{paragraph}\noi 
Let $D(R)$ denote the derived category of the $R$-module category $R$-Mod. A complex $M \in D(R)$ is called a {\it derived $I$-torsion complex} if the
canonical morphism ${\bf R}\Gamma_I(M)\rightarrow M$ is an isomorphism. The complex $M$ is called
a {\it derived $I$-adically complete complex} if the canonical morphism $M\rightarrow {\bf L}\Lambda_I(M)$ is an isomorphism. Denote by $D(R)_{I\text{-tor}}$ and $D(R)_{I\text{-com}}$ the full
subcategories of $D(R)$ consisting of derived $I$-torsion complexes and derived $I$-adically complete complexes, respectively. These are triangulated subcategories.  One version of the MGM Equivalence is Theorem \ref{q} which appears as Theorem 7.11 in  \cite{Yekutieli}. Recently, its version in the noncommutative setting was given, see  \cite{Vyas2}.  
\end{paragraph}

    \begin{thm}\label{q}  
    Let $R$ be a ring, and let $I$ be a weakly proregular ideal in it.

\begin{enumerate}
 \item  For any $M\in {\bf D}(R)$,  
 ${\bf R}\Gamma_I(M)\in {\bf D}(R)_{I\text{-tor}}~ \text{ and} ~~{\bf L}\Lambda_I(M)\in {\bf D}(R)_{I\text{-com}}.$
 
\item  The functor $$ {\bf R}\Gamma_I~: ~ {\bf D}(R)_{I\text{-com}}\rightarrow {\bf D}(R)_{I\text{-tor}}$$ is an  equivalence, with quasi-inverse ${\bf L}\Lambda_I$.
\end{enumerate}

\end{thm}

\begin{paragraph}\noi 
 Our second main result is Theorem \ref{02} which realises the MGM Equivalence in the setting of $I$-reduced and $I$-coreduced $R$-modules. An $R$-module is $I$-torsion (resp. $I$-complete) if $\Gamma_I(M)\cong M$ (resp. $\Lambda_I(M)\cong M$). The full subcategory of $R$-Mod consisting of $I$-torsion (resp. $I$-complete) $R$-modules is denoted by $(R\text{-Mod})_{I\text{-tor}}$ (resp. $(R\text{-Mod})_{I\text{-com}}$).
\end{paragraph}

\begin{thm}{\rm[{\bf The MGM Equivalence in $R$-Mod}]}\label{02} Let $I$ be any ideal of a ring $R$, 
   
   \begin{enumerate}
    \item For any $M\in (R\text{-Mod})_{I\text{-red}}$, $$\Gamma_I(M)\in  (R\text{-Mod})_{I\text{-com}}\cap (R\text{-Mod})_{I\text{-cor}}=:\mathfrak{E}.$$

     \item For any $M\in (R\text{-Mod})_{I\text{-cor}}$, $$\Lambda_I(M)\in  (R\text{-Mod})_{I\text{-tor}}\cap (R\text{-Mod})_{I\text{-red}}=:\mathfrak{D}.$$
  
\item The functor $\Gamma_I : (R\text{-Mod})_{I\text{-red}}\rightarrow (R\text{-Mod})_{I\text{-cor}}$ restricted to $\mathfrak{D}$ is an equivalence between $\mathfrak{D}$ and $\mathfrak{E}$ with quasi-inverse $\Lambda_I$.    
     \end{enumerate}
   \end{thm}

\begin{paragraph}\noi

Reduced modules were introduced by Lee and Zhou in \cite{Lee} and later studied in \cite{Rege}. In both papers, the method of study was mainly element-wise.  It was however observed in \cite{David} that  ``reduced modules'' is a categorical property. Therefore, it is amenable to study by use of category theory. The use of this method to study reduced modules first appeared in \cite{Kyom}. In this paper, we  continue with this approach; it allows us to benefit from the powerful machinery of category theory. The key idea that makes all that we do possible is the fact that the functor $\Gamma_I$ when restricted to $I$-reduced $R$-modules is representable, i.e., for all $I$-reduced $R$-modules $M$, $\Gamma_I(M)\cong \text{Hom}_R(R/I, M)$. Dually, for any $I$-coreduced $R$-module $M$, $\Lambda_I(M)$ is naturally isomorphic to $R/I\otimes M$.  

\end{paragraph}

  \begin{paragraph}\noi 
   The paper is organised as follows. It consists of five sections. In Section 2, we give the basic properties of $I$-reduced and $I$-coreduced $R$-modules (given an ideal $I$ of $R$) necessary for the subsequent sections. In Sections 3 and 4, we prove the GM Duality and the MGM Equivalence respectively. Weakly proregular ideals have been hitherto known to be the most general condition under which the GM Duality and the MGM Equivalence hold (although in the derived category setting). Since in  this paper we have another condition of $I$-reduced and $I$-coreduced $R$-modules serving the same purpose, part of Section 4 aims at comparing the two conditions. In general, none of the conditions implies the other. The last section; Section 5, explores more ways (beyond those in \cite{Kyom} for $I$-reduced $R$-modules and  in \cite{Kyom2} for $I$-coreduced $R$-modules) of realising the functors $\Gamma_I$ and $\Lambda_I$ as representable functors. Furthermore, we demonstrate that doing so is important by utilising the Yoneda Lemma to compute natural transformations.
  \end{paragraph}

\section{Basic properties}

\begin{defn}\rm Let $I$ be an ideal of a ring $R$, $M$  be an $R$-module and $f_a$ be the endomorphism of $M$ given by $m\mapsto am$ for $a\in R$. $M$ is 
 \begin{enumerate}
  \item {\it $I$-reduced}  if for every   $a\in I$, $ \text{Ker}~f_a \cong \text{Ker}~f_a^2;$
   
  \item  {\it $I$-coreduced} if for every  $a\in I$, $ \text{Coker}~f_a \cong \text{Coker}~f_a^2.$   
   \end{enumerate}

 \end{defn}

 \begin{paragraph}\noi
 An $R$-module $M$ is {\it reduced} (resp. {\it coreduced}) if for every ideal $I$ of $R$, $M$ is $I$-reduced (resp. $I$-coreduced). A ring $R$ is reduced (resp. coreduced) if and only if it is reduced (resp. coreduced) as an $R$-module.  Since a ring is reduced if and only if it has no nonzero nilpotent ideals and it is coreduced if and only if all its ideals are idempotent, it follows that idempotent ideals dualise having no nonzero nilpotent ideals.   
  \end{paragraph}

 \begin{paragraph}\noi 
  For any ideal $I$ of a ring $R$, the class of $I$-reduced $R$-modules is quite large. It contains reduced $R$-modules which also contain the well studied class of prime $R$-modules. Dually, the $I$-coreduced $R$-modules contain the class of coreduced $R$-modules which in turn contain the second (also called coprime) $R$-modules.
 \end{paragraph}

\begin{propn}\label{pr}
 For any $R$-module $M$ and an ideal $I$ of $R$, the following statements are equivalent:
 
\begin{enumerate}
 \item $M$ is $I$-reduced,
 \item  $(0:_MI)=(0:_MI^2)$,
 \item $\text{Hom}_R(R/I, M)\cong \text{Hom}_R(R/I^2, M)$,
 \item $\Gamma_I(M)\cong \text{Hom}_R(R/I, M)$, 
 \item $I\Gamma_I(M)=0$.
 \end{enumerate}

\end{propn}

\begin{prf}
 \begin{itemize}
  \item[$1\Rightarrow 2$] Since in general $(0:_MI)\subseteq (0:_MI^2)$, let $m\in (0:_MI^2)$. Then $I^2m=0$ and $a^2m=0$ for all $a\in I$. So, $m\in (0:_Ma^2)$ for all $a\in I$.   $\text{Ker} f_a^2=(0:_Ma^2)$ and by hypothesis, $\text{Ker}~f_a=\text{Ker}~f_a^2$. So, we get $m\in (0:_Ma)$ for all $a\in I$ and hence $m\in (0:_MI)$.
  
  \item[$2\Rightarrow 3$] This is immediate since $\text{Hom}_R(R/I, M)\cong (0:_MI)$.
  \item[$3\Rightarrow 4$] By definition, $\Gamma_I(M):= \lim\limits_{\rightarrow}\text{Hom}_R(R/I^k, M)$. It follows by statement (3) that $\text{Hom}_R(R/I^k, M)\cong \text{Hom}_R(R/I, M)$ for all $k\in \Z^+$. So, $\Gamma_I(M)\cong \text{Hom}_R(R/I, M)$.
  \item[$4\Rightarrow 5$] If $\Gamma_I(M)\cong \text{Hom}_R(R/I, M)$, then $I\Gamma_I(M)\cong I\text{Hom}_R(R/I, M)\cong \text{Hom}_R(0, M)\cong 0$.
  \item[$5\Rightarrow 1$] Since for any $a\in I$, $\text{Ker}~f_a\subseteq~\text{Ker}~f_a^2$, we prove the reverse inclusion. Let $m\in \text{Ker}~f_a^2=(0:_Ma^2)$ for all $a\in I$. $a^2m=0$ and $(a)^2m=0$ for all $a\in I$, where $(a)$ is the ideal of $R$ generated by $a$. This implies that $m\in\Gamma_{(a)}(M)$ for all $a\in I$ and $m\in \bigcap\limits_{a\in I}\Gamma_{(a)}(M)=\Gamma_I(M)$. Therefore, $am\in Im\subseteq I\Gamma_I(M)=0$ for all $a\in I$ and hence $m\in (0:_Ma)=\text{Ker}~f_a$ for all $a\in I$.
 \end{itemize}

\end{prf}

\begin{propn}\label{pc}
 For any $R$-module $M$ and an ideal $I$ of $R$, the following statements are equivalent:
 
\begin{enumerate}
 \item $M$ is $I$-coreduced,
 \item $IM=I^2M$,
 \item $R/I\otimes M\cong R/I^2\otimes M$,
 \item $\Lambda_I(M)\cong R/I\otimes M$, 
 \item $I\Lambda_I(M)=0$.
  
\end{enumerate}
\end{propn}

\begin{prf}
 \begin{itemize}
  \item[$1\Rightarrow 2$]Suppose that $M/a^2M\cong M/aM$ for all $a\in I$. So,  $a^2M=aM$ for all $a\in I$. To see this, note that the natural epimorphisms $M\rightarrow M/a^2M$ and $M\rightarrow M/aM$ have kernels $a^2M$ and $aM$ respectively. The isomorphism $M/a^2M\cong M/aM$ implies that the submodules $a^2M$ and $aM$ of $M$ should coincide for all $a\in I$. Therefore,  $I^2M=IM$.
  \item[$2\Rightarrow 3$] $R/I\otimes M\cong M/IM=M/I^2M \cong R/I^2\otimes M.$
  \item[$3\Rightarrow 4$] $\Lambda_I(M):=\lim\limits_{\leftarrow}(R/I^k\otimes M)\cong \lim\limits_{\leftarrow}(R/I\otimes M)=R/I\otimes M$.
  \item[$4\Rightarrow 5$] Given statement (4), we have $I\Lambda_I(M)\cong I(R/I\otimes M)\cong (0\otimes M)\cong 0$.
  \item[$5\Rightarrow 1$]$0\cong I\Lambda_I(M)=I\lim\limits_{\leftarrow}(M/I^kM)$ implies that $IM=I^kM$ for all $k\in \Z^+$. Since $a^2M\subseteq aM$ for all $a\in I$, let $m\in aM$ for all $a\in I$. Then $m\in IM=I^kM$ for all $k\in \Z^+$. So, $m\in a^2M$ for all $a\in I$.
  This establishes the equality $a^2M=aM$ for all $a\in I$.  Thus $M/aM=M/a^2M$ for all $a\in I$.
 \end{itemize}

\end{prf}
 
\begin{rem}\rm 
 Proposition \ref{pr} holds when in the place of the ideal $I$ of $R$ one has an element $a\in R$ (or the principal ideal $aR$) in which case $M$ is called $a$-reduced, see \cite[Proposition 2.2]{Kyom}. Proposition \ref{pc} already appears in \cite{Kyom2}. For completeness, we have given the proof.
\end{rem}

 \begin{propn}\label{cisr}
  Every coreduced ring is reduced.
 \end{propn}
 
 \begin{prf}
  If a ring $R$ is coreduced, then every ideal $I$ of $R$ is idempotent since $I^2R=IR$. If $I^2=0$, then $I=0$ and $R$ is reduced.
 \end{prf}

 \begin{paragraph}\noi 
 The converse of Proposition \ref{cisr} is not true in general. The ring of integers is reduced but it is not coreduced. A coreduced module also need not be reduced. 
 The $\Z$-module $\Q/\Z$ is coreduced but it is not reduced.
 \end{paragraph}

\begin{propn}\label{P1}
Let $R$ be  a ring. 

\begin{enumerate}
 \item  For any $R$-module $N$ and a coreduced (resp. $I$-coreduced) $R$-module $M$, the $R$-module $\text{Hom}_R(M,N)$ is reduced (resp. $I$-reduced).
 
 \item   Let $N$ be an injective cogenerator of  $R$-Mod. If $\text{Hom}_R(M,N)$ is a reduced (resp. $I$-reduced) $R$-module for some $M\in R$-Mod, then $M$ is a coreduced (resp. $I$-coreduced) $R$-module.
 
 \item \label{cog}
Let $N$ be an injective cogenerator $R$-module. The $R$-module $M$ is coreduced (resp. $I$-coreduced) if and only if $\text{Hom}_R(M,N)$ is a reduced (resp. $I$-reduced) $R$-module.
\end{enumerate}
 \end{propn}

\begin{prf}  We prove the $I$-reduced and $I$-coreduced cases. The reduced and coreduced versions follow immediately. 
\begin{enumerate} 
\item  Let $I$ be an ideal of $R$ and $M$ an $I$-coreduced $R$-module. By the Hom-Tensor duality and Proposition \ref{pc}, we have  $\text{Hom}_R(R/I, \text{Hom}_R(M,N))\cong \text{Hom}_R(R/I\otimes M, N)$
 $\cong \text{Hom
 }_R(R/I^2\otimes M, N)$ $\cong \text{Hom}_R(R/I^2, \text{Hom}_R(M, N))$. By Proposition \ref{pr}, $\text{Hom}_R(M,N)$ is $I$-reduced.

\item   Let $I$ be an ideal of $R$ and $\text{Hom}_R(M, N)$ be an $I$-reduced $R$-module.  By the Hom-Tensor duality and  Proposition \ref{pr},
 $$  \text{Hom}_R(R/I\otimes M, N)
 \cong \text{Hom}_R(R/I, \text{Hom}_R(M,N)) $$
 $$\cong \text{Hom}_R(R/I^2, \text{Hom}_R(M, N)) \cong\text{Hom}_R(R/I^2\otimes M, N).$$ Since $N$ is an injective cogenerator, the functor $\text{Hom}_R(-,N)$  reflects isomorphisms. So, $R/I\otimes M\cong R/I^2\otimes M$ and by Proposition \ref{pc}, $M$ is $I$-coreduced.
 
 \item This is immediate from 1) and 2) above. 
\end{enumerate}
\end{prf}

\begin{propn}\label{P3}
 Let $M$ be an $S$-$R$-bimodule and $N$ a left $R$-module. If $J$ is an ideal of $S$ and $_SM$ is a $J$-coreduced (resp. coreduced) $S$-module, then so is the $S$-module $M\otimes N$.
\end{propn}

\begin{prf}
 If $J$ is an ideal of $S$, then $S/J \otimes_S (M\otimes _RN)\cong(S/J\otimes_SM)\otimes_RN\cong(S/{J^2}\otimes_SM)\otimes_RN\cong S/{J^2}\otimes_S(M\otimes_RN)$.
\end{prf}

 \begin{cor}Let $L$ be an $R$-module.
  If $\text{Hom}_R(L,R)$ is a coreduced $R$-module, then so is the $R$-module $\text{Hom}_R(L, M)$ for any $M\in \text{R-Mod}$.
 \end{cor}
 
 \begin{prf}
  Follows from Proposition \ref{P3} and the fact that $\text{Hom}_R(L,M)\cong \text{Hom}_R(L,R)\otimes M$.
 \end{prf}

  \begin{paragraph}\noi   
  Let $I$ be an ideal of a ring $R$. We denote by $(R$-Mod$)_{I\text{-cor}}$ (resp. $(R$-Mod$)_{I\text{-red}}$) the subcategory of  $R\text{-Mod}$ consisting of  $I$-coreduced (resp. $I$-reduced) $R$-modules.
   
  \end{paragraph}

  \begin{propn}\label{cr}
   For any ideal $I$ of a ring $R$, we have:
   
   \begin{enumerate}
    \item $R/I\otimes - $ and $\text{Hom}_R(R/I, -)$ are  idempotent functors from $R$-Mod to
    $$(R\text{-Mod})_{I\text{-cor}} \cap (R\text{-Mod})_{I\text{-red}}.$$
    
    \item For any $R$-module $M$, 
    
    \begin{equation}
     R/I\otimes \text{Hom}_R(R/I, M)\cong \text{Hom}_R(R/I, M).
    \end{equation}
    
    \begin{equation}
     \text{Hom}_R(R/I, R/I\otimes M)\cong R/I\otimes M.     
    \end{equation}
    
    \item For any $R$-module $M$, the $R$-modules $\text{Hom}_R(R/I, M)$ and $R/I\otimes M$ are both $I$-torsion and $I$-complete.
    
    \item The set $\mathfrak{A}:=\{\text{Hom}_R(R/I, -), R/I\otimes -\}$ forms a noncommutative semigroup where the operation is composition of functors.
   \end{enumerate}
  \end{propn}

  \begin{prf}
  
  \begin{enumerate}
   \item 
    Idempotency holds because
  $$R/I\otimes (R/I\otimes M)\cong (R/I\otimes R/I)\otimes M\cong R/I\otimes M$$
  and 
  $$\text{Hom}_R(R/I, \text{Hom}_R(R/I,M))\cong \text{Hom}_R(R/I\otimes R/I, M)\cong \text{Hom}_R(R/I,M).$$

  For any $R$-module $M$, $R/I\otimes M \cong M/IM$ and $\text{Hom}_R(R/I, M)\cong (0:_MI)$. Also, $I(M/IM)=0$  and $I(0:_MI)=0$. So, the $R$-modules $M/IM$  and $(0:_MI)$ are  $I$-coreduced.   
   It is also easy to see that 
   $$(0:_{(0:_MI)}I)= (0:_{(0:_MI)}I^2)=(0:_MI)$$ and $$(\bar{0}:_{M/IM}I)= (\bar{0}:_{M/IM}I^2)=M/IM.$$ This shows that the $R$-modules $M/IM$  and $(0:_MI)$ are  $I$-reduced.
   
   \item $$R/I\otimes \text{Hom}_R(R/I,M)\cong \frac{\text{Hom}_R(R/I, M)}{I\text{Hom}_R(R/I,M)}\cong \frac{\text{Hom}_R(R/I, M)}{0}\cong \text{Hom}_R(R/I,M).$$
   
   $$\text{Hom}_R(R/I, R/I\otimes M)\cong \text{Hom}_R(R/I, M/IM)\cong (\bar{0}:_{M/IM}I)=M/IM \cong R/I\otimes M.$$

   \item The following maps always hold\footnote{$\hookrightarrow$ denotes a monomorphism and $\twoheadrightarrow$ denotes an epimorphism.} 
   
   \begin{equation}\label{1}
   \text{Hom}_R(R/I, \text{Hom}_R(R/I,M))\hookrightarrow \Gamma_I(\text{Hom}_R(R/I, M))\hookrightarrow \text{Hom}_R(R/I, M)    
   \end{equation}
   
   \begin{equation}\label{2}
   \text{Hom}_R(R/I, R/I\otimes M)\hookrightarrow \Gamma_I(R/I\otimes M)\hookrightarrow R/I\otimes M   
   \end{equation}
   
   \begin{equation}\label{3}
    \text{Hom}_R(R/I, M) \twoheadrightarrow \Lambda_I(\text{Hom}_R(R/I, M) ) \twoheadrightarrow R/I\otimes \text{Hom}_R(R/I, M) 
   \end{equation}
   
   \begin{equation}\label{4}
    R/I\otimes M \twoheadrightarrow \Lambda_I(R/I\otimes M ) \twoheadrightarrow R/I\otimes (R/I\otimes M)
   \end{equation}
   
   Idempotency of the functors $\text{Hom}_R(R/I, -)$ and $R/I\otimes -$ shows that the maps in (\ref{1}) and (\ref{4}) are all isomorphisms in which case, $\text{Hom}_R(R/I, M)$ and $R/I\otimes M$ become $I$-torsion and $I$-complete respectively.  Invariance of $R/I\otimes M$ and $\text{Hom}_R(R/I, M)$ under the functor $\text{Hom}_R(R/I, -)$ and $R/I\otimes -$ respectively shows that the morphisms in (\ref{2}) and (\ref{3}) are all isomorphisms. This shows that $R/I\otimes M$ and $\text{Hom}_R(R/I, M)$ are $I$-torsion and $I$-complete respectively.

   \item From 1) and 2) above, we get Figure \ref{fg} which shows that the set $\mathfrak{A}$ is a noncommutative semigroup.
   
   \begin{figure}[h!]
  \begin{center}
     \begin{tikzpicture}   
 \draw[thick, -]  (0,0) -- (8, 0);
 \draw[thick, -]  (0,-1) -- (8, -1);
 \draw[thick, -]  (3, 1) -- (3, -2);
 \draw[thick, -]  (6, 1) -- (6, -2);
  \node (A) at (4.5, 0.5) {$\text{Hom}_R(R/I, -)$}; 
  \node (A) at (7, 0.5) {$R/I\otimes -$}; 
  \node (A) at (4.5, -0.5) {$\text{Hom}_R(R/I, -)$}; 
  \node (A) at (7, -0.5) {$R/I\otimes -$}; 
    \node (A) at (4.5, -1.5) {$\text{Hom}_R(R/I, -)$}; 
  \node (A) at (7, -1.5) {$R/I\otimes -$};
   \node (A) at (1.5, -0.5) {$\text{Hom}_R(R/I, -)$}; 
  \node (A) at (1.5, -1.5) {$R/I\otimes -$};     
   \end{tikzpicture}    
  \end{center}
   \caption{Multiplication table} \label{fg}
   \end{figure} 
   
   \end{enumerate}
  \end{prf}

  \begin{cor}
   Every $R$-module $M$ has  a submodule and a quotient module which are both $I$-reduced and  $I$-coreduced as $R$-modules.
  \end{cor}
  
  \begin{prf}
   Immediate from Proposition \ref{cr}(1). They are $(0:_MI)$ and $M/IM$ respectively.
  \end{prf}

  \begin{cor}
   If $M$  is $I$-reduced (resp. $I$-coreduced), then $$\text{Hom}_R(R/I, \Gamma_I(M))\cong \text{Hom}_R(R/I, M)\cong \Gamma_I(M)$$  and 
   $$ (\text{resp.}~ R/I\otimes \Lambda_I(M) \cong R/I\otimes M \cong \Lambda_I(M)).$$
  \end{cor}

  \begin{prf}
   By Proposition \ref{pr}, $\Gamma_I(M)\cong \text{Hom}_R(R/I, M)$. So, $$\text{Hom}_R(R/I, \Gamma_I(M))\cong \text{Hom}_R(R/I, \text{Hom}_R(R/I, M))\cong \text{Hom}_R(R/I, M)\cong \Gamma_I(M).$$ Also by Proposition \ref{pc}, $\Lambda_I(M)\cong R/I\otimes M$ and therefore $R/I\otimes \Lambda_I(M) \cong R/I\otimes R/I\otimes M \cong R/I\otimes M \cong \Lambda_I(M)$.
  \end{prf}

  \begin{propn}\label{limits} Let $I$ be an ideal of a ring $R$. An inverse limit of $I$-reduced (resp. reduced) $R$-modules is an $I$-reduced (resp. reduced) $R$-module and a direct limit of $I$-coreduced (resp. coreduced) $R$-modules is an $I$-coreduced (resp. coreduced) $R$-module.

  \end{propn}

  \begin{prf} We prove the cases for $I$-reduced and $I$-coreduced. The other cases become immediate.
   It is well known that the functor $\text{Hom}_R(R/I, -)$ (resp. $R/I\otimes -)$ preserves inverse limits (resp. direct limits). A property possessed by right adjoint (resp. left adjoint) functors, see for instance \cite[V. 5]{Mac}. So, if each $R$-module $M_k$ is $I$-reduced, then  $\text{Hom}_R(R/I, \lim\limits_{\leftarrow}M_k)\cong \lim\limits_{\leftarrow}\text{Hom}_R(R/I, M_k)\cong \lim\limits_{\leftarrow}\text{Hom}_R(R/I^2, M_k)\cong \text{Hom}_R(R/I^2, \lim\limits_{\leftarrow}M_k)$ which shows that
   $\lim\limits_{\leftarrow}M_k$ is an $I$-reduced $R$-module. Dually, suppose that each $R$-module $M_k$ is $I$-coreduced. $R/I\otimes \lim\limits_{\rightarrow}M_k\cong \lim\limits_{\rightarrow}(R/I\otimes M_k)\cong \lim\limits_{\rightarrow}(R/I^2\otimes M_k)=R/I^2\otimes \lim\limits_{\rightarrow}M_k$ so that $\lim\limits_{\rightarrow}M_k$ is an $I$-coreduced $R$-module.
  \end{prf}

\begin{cor}\label{cr1}
 Let $I$ be an ideal of $R$ and $M$ be an $R$-module. $M$ is $I$-reduced if and only if so  is the $R$-submodule $\Gamma_I(M)$ of $M$.
\end{cor}

\begin{prf}
 $M$ is $I$-reduced if and only if $I\Gamma_I(M)=0$ (Proposition \ref{pr}). Idempotency of $\Gamma_I$ implies that  $I\Gamma_I(M)=0$ if and only if  $I\Gamma_I(\Gamma_I(M))=0$. This is the case if and only if $\Gamma_I(M)$ is $I$-reduced (Proposition \ref{pr}).
\end{prf}

\begin{paragraph}\noi 
 By  a similar proof, we obtain:
\end{paragraph}

\begin{cor}\label{cr2}
 Let $I$ be an ideal of $R$,  $M$ be an $R$-module and $\Lambda_I$ be an idempotent functor, then $M$ is $I$-coreduced if and only if so  is  $\Lambda_I(M)$.
\end{cor}

\section{Greenless-May Duality}

 \begin{paragraph}\noi 
 For an arbitrary ideal $I$ of a ring $R$ and $R$-modules $M$ and $N$ $$\text{Hom}_R(\Lambda_{I}(M), N)\not\cong \text{Hom}_{R}(M, \Gamma_{I}(N)).$$
 So, the functors $\Gamma_I$ and $\Lambda_I$ are not in general adjoint to each other.  However, in the setting of derived categories, we have Theorem \ref{T3} which is called the Greenless-May Duality (GM Duality for short). It was first proved by  Alonso Tarrio, Jeremias Lopez and Lipman \cite{Alonso} but it also appears in \cite[Theorem  7.1.2]{Yekutieli}.
\end{paragraph}
 
\begin{thm}\label{T3}\rm [GM-Duality in $\mathbf{D}(R)$] Let $I$ be a weakly proregular ideal of a ring $R$ and $M, N\in \mathbf{D}(R)$. Then there is a natural isomorphism in $\mathbf{D}(R)$ given by $$\mathbf{R}\text{Hom}_R(\mathbf{R}\Gamma_I(M), N)\cong \mathbf{R}\text{Hom}_R(M, \mathbf{L}\Lambda_I(N)).$$ 
 \end{thm}
\begin{paragraph}\noi  
This theorem implies that local cohomology is derived left adjoint to local homology.  Let $I$ be an ideal of a ring $R$ and let  $(R\text{-Mod})_{I\text{-red}}$  (resp. $(R\text{-Mod})_{I\text{-cor}}$) denote a full subcategory of $R$-Mod consisting of all $I$-reduced (resp. $I$-coreduced) $R$-modules. 
\end{paragraph}

\begin{thm}\label{TGM}{\rm {\bf [GM Duality in $R$-Mod]}}
 For any ideal $I$ of a ring $R$, 
 
 \begin{enumerate}
  \item  The functor 
 $$\Gamma_I: (R\text{-Mod})_{I\text{-red}} \rightarrow (R\text{-Mod})_{I\text{-cor}}$$ is idempotent and for any $M\in (R\text{-Mod})_{I\text{-red}}$, $\Gamma_I(M)\cong \text{Hom}_R(R/I, M)$.
 
 \item   The functor $$\Lambda_I: (R\text{-Mod})_{I\text{-cor}} \rightarrow (R\text{-Mod})_{I\text{-red}}$$ is idempotent and  for any  $M\in (R\text{-Mod})_{I\text{-cor}}$, $\Lambda_I(M)\cong R/I\otimes M$.
 \item  For any $N\in (R\text{-Mod})_{I\text{-red}}$ and $M\in (R\text{-Mod})_{I\text{-cor}}$, 
 
 $$\text{Hom}_R(\Lambda_I(M), N)\cong \text{Hom}_R(M, \Gamma_I(N)).$$  

 \end{enumerate}

\end{thm}

\begin{prf}

\begin{enumerate}
\item  From Proposition \ref{pr}, $\Gamma_I(M)$ is naturally isomorphic to $\text{Hom}_R(R/I, M)$ for any $I$-reduced $R$-module $M$. Idempotency is due to Proposition \ref{cr}(1) and the fact that $\Gamma_I(M)\cong\text{Hom}_R(R/I, M)$ for all $I$-reduced $R$-modules $M$.

\item   Follows from Proposition \ref{pc}(2). Idempotency is due to Proposition \ref{cr}(1) and the fact that $\Lambda_I(M)\cong R/I\otimes M$ for all $I$-coreduced $R$-modules $M$.
 
 \item Consider the functor $\Gamma_I: (R\text{-Mod})_{I\text{-red}} \rightarrow (R\text{-Mod})_{I\text{-cor}}$. For any module $M\in (R\text{-Mod})_{I\text{-red}} $, $\Gamma_I(M)\cong \text{Hom}_R(R/I, M)$. However, the functor $R/I\otimes -$ is left-adjoint to $\text{Hom}_R(R/I, -)$. By uniqueness of adjoints, the functor
$ \Lambda_I: (R\text{-Mod})_{I\text{-cor}} \rightarrow (R\text{-Mod})_{I\text{-red}}$ which has the property that for all $M\in 
     (R\text{-Mod})_{I\text{-cor}}$, $\Lambda_I(M)\cong R/I\otimes M$ is the left adjoint of $\Gamma_I$.      
\end{enumerate}

\end{prf}

\begin{cor}\label{ccr}
 If $M\in (R\text{-Mod})_{I\text{-red}}\cap (R\text{-Mod})_{I\text{-cor}}$, then $\Gamma_I(M)=0$ if and only if $\Lambda_I(M)=0$.
\end{cor}

\begin{prf}
 By the adjunction in Theorem \ref{TGM}, if $M=N$, we get 
 $$\text{Hom}_R(\Lambda_I(M), M)\cong \text{Hom}_R(M, \Gamma_I(M)).$$
 If $\Lambda_I(M)=0$, $\text{Hom}_R(M, \Gamma_I(M))=0$, and every $R$-homomorphism $f:M\rightarrow \Gamma_I(M)$ is a zero homomorphism. Applying $\Gamma_I$ gives $\Gamma_I(f):\Gamma_I(M)\rightarrow\Gamma_I(\Gamma_I(M))$ which is zero. However, $\Gamma_I$ is idempotent. So $\Gamma_I(f)$ is an isomorphism and therefore  $\Gamma_I(M)=0$. Conversely, suppose that $\Gamma_I(M)=0$. Then $\text{Hom}_{\mathfrak{R}}(\Lambda_I(M), M)=0$ and therefore every $R$-homomorphism $g: \Lambda_I(M)\rightarrow M$ is zero.
 It follows that $\Lambda_I(g): \Lambda_I(\Lambda_I(M))\rightarrow \Lambda_I(M)$ is also a zero $R$-homomorphism.  By idempotency of $\Lambda_I$ on $I$-coreduced $R$-modules $\Lambda_I(g)$ is an isomorphism and hence $\Lambda_I(M)=0$.
 
\end{prf}

\begin{exam}\rm We have the following examples.
 \begin{enumerate}
  \item If $I^2=I$, then $R\text{-Mod}= (R\text{-Mod})_{I\text{-red}} = (R\text{-Mod})_{I\text{-cor}}$. 
  
  \item For any simple $R$-module $M$, $M\in (R\text{-Mod})_{I\text{-red}}\cap (R\text{-Mod})_{I\text{-cor}}$.
  
  \item If $R$ is an Artinian ring with an ideal $I$, then there exists a positive integer $k$ such that every $R$-module is both $I^k$-reduced and $I^k$-coreduced, i.e., $R\text{-Mod}= (R\text{-Mod})_{I^k\text{-red}}\cap (R\text{-Mod})_{I^k\text{-cor}}$. 
 \end{enumerate}
The first two examples are easy to see. Suppose that $R$ is Artinian. There exists some positive integer  $k$ such that $I^k=I^{k+t}$ for all positive integers $t$. So, $(0:_MI^k)=(0:_MI^{k+t})$ and $M/I^kM=M/I^{k+t}M$ for all positive integers $t$. Accordingly, 
$\Gamma_{I^k}(M)\cong (0:_MI^k)$ and $\Lambda_{I^k}(M)\cong M/I^kM$. So, $M$ is both $I^k$-reduced and $I^k$-coreduced.
\end{exam}

\begin{cor}Let $I$ be an ideal of a ring $R$, 
\begin{enumerate}
 \item For any $M\in (R\text{-Mod})_{I\text{-red}}$, $\Lambda_I(\Gamma_I(M))=\Gamma_I(M)$.
 
 \item For any $M\in (R\text{-Mod})_{I\text{-cor}}$, $\Gamma_I(\Lambda_I(M))=\Lambda_I(M)$.
 \end{enumerate}
 
 \begin{prf}
 \begin{enumerate}
  \item 
  By Theorem \ref{TGM}, for any $M\in (R\text{-Mod})_{I\text{-red}}$, $\Gamma_I(M)\in (R\text{-Mod})_{I\text{-cor}}$ and $\Gamma_I(M)\cong \text{Hom}_R(R/I,M)$. So, $$\Lambda_I(\Gamma_I(M))\cong R/I\otimes \text{Hom}_R(R/I, M)\cong \text{Hom}_R(R/I, M) \cong \Gamma_I(M).$$
  
  \item By Theorem \ref{TGM}, for any $M\in (R\text{-Mod})_{I\text{-cor}}$, $\Lambda_I(M)\in (R\text{-Mod})_{I\text{-red}}$ and $\Lambda_I(M)\cong R/I\otimes M$. It follows that 
  $\Gamma_I(\Lambda_I(M)) \cong \text{Hom}_R(R/I, R/I\otimes M)\cong R/I\otimes M \cong \Lambda_I(M)$.
  \end{enumerate}
 \end{prf}
\end{cor}

    \section{MGM Equivalence}
    
    \begin{paragraph}\noi 
    Let $R$ be a ring and ${\bf D}(R)$ denote the derived category of the abelian category $R$-Mod.  
    A complex $M \in {\bf D}(R)$ is called  
    {\it derived $I$-torsion} \cite[Definition 3.11]{Yekutieli} (resp.  {\it  derived $I$-adically complete} \cite[Definition 3.8]{Yekutieli}) if the morphism    
    
    $$ \sigma_M^R : {\bf R}\Gamma_I(M) \rightarrow  M~~~ (\text{resp.}~\tau_M^L : M\rightarrow {\bf L}\Lambda_I(M) )$$ is an isomorphism, where ${\bf R}\Gamma_I$ and ${\bf L}\Lambda_I$ denote the right derived and left derived functors of $\Gamma_I$ and $\Lambda_I$ respectively.
The full subcategory of ${\bf D}(R)$ consisting of derived $I$-torsion (resp.  derived $I$-adically complete) complexes is denoted by 
${\bf D}(R)_{I\text{-tor}}$ (resp. ${\bf D}(R)_{I\text{-com}}$).
\end{paragraph}

    \begin{thm}\rm[{\bf MGM Equivalence \cite[Theorem 7.11]{Yekutieli}]}    
    Let $R$ be a ring, and let $I$ be a weakly proregular ideal in it.

\begin{enumerate}
 \item  For any $M\in {\bf D}(R)$,  
 ${\bf R}\Gamma_I(M)\in {\bf D}(R)_{I\text{-tor}}~ \text{ and} ~~{\bf L}\Lambda_I(M)\in {\bf D}(R)_{I\text{-com}}.$
 
\item  The functor $$ {\bf R}\Gamma_I~: ~ {\bf D}(R)_{I\text{-com}}\rightarrow {\bf D}(R)_{I\text{-tor}}$$ is an  equivalence, with quasi-inverse ${\bf L}\Lambda_I$.
\end{enumerate}

\end{thm}

\begin{lem}\label{rediscom}
 If $I$ is an ideal of a ring $R$ and $M$ an $I$-reduced (resp. $I$-coreduced) $R$-module, then $\Gamma_I(M)$ (resp. $\Lambda_I(M)$) is an $I$-complete (resp. $I$-torsion) $R$-module.
   \end{lem}
   
   \begin{prf}
    If $M$ is $I$-reduced, $\Gamma_I(M)\cong \text{Hom}_R(R/I, M)$ which is both an $I$-reduced and $I$-coreduced $R$-module by Proposition \ref{cr}. So, $\Lambda_I(\Gamma_I(M))\cong R/I\otimes \text{Hom}_R(R/I, M)\cong \text{Hom}_R(R/I, M)\cong \Gamma_I(M)$. This proves that $\Gamma_I(M)$ is $I$-complete. If $M$ is $I$-coreduced, $\Lambda_I(M)\cong R/I \otimes M $ which is also both an $I$-reduced and $I$-coreduced $R$-module by Proposition \ref{cr}. So, $\Gamma_I(\Lambda_I(M))\cong  \text{Hom}_R(R/I, \Lambda_I(M))\cong \text{Hom}_R(R/I, R/I\otimes M)\cong R/I\otimes M\cong\Lambda_I(M)$. This proves that $\Lambda_I(M)$ is $I$-torsion. 
   \end{prf}

   \begin{thm}{\rm[{\bf The MGM Equivalence in $R$-Mod}]}\label{MGM} Let $I$ be any ideal of a ring $R$, 
   
   \begin{enumerate}
    \item For any $M\in (R\text{-Mod})_{I\text{-red}}$, $$\Gamma_I(M)\in  (R\text{-Mod})_{I\text{-com}}\cap (R\text{-Mod})_{I\text{-cor}}=:\mathfrak{E}.$$

     \item For any $M\in (R\text{-Mod})_{I\text{-cor}}$, $$\Lambda_I(M)\in  (R\text{-Mod})_{I\text{-tor}}\cap (R\text{-Mod})_{I\text{-red}}=:\mathfrak{D}.$$
  
\item The functor $\Gamma_I : (R\text{-Mod})_{I\text{-red}}\rightarrow (R\text{-Mod})_{I\text{-cor}}$ restricted to $\mathfrak{D}$ is an equivalence between $\mathfrak{D}$ and $\mathfrak{E}$ with quasi-inverse $\Lambda_I$.

\item If $I$ is an idempotent ideal, then the subcategory $(R\text{-Mod})_{I\text{-com}}$ is equivalent to $(R\text{-Mod})_{I\text{-tor}}$.
    
     \end{enumerate}
   \end{thm}

\begin{prf}
\begin{enumerate}
 \item By Theorem \ref{TGM}, $M\in( R\text{-Mod})_{I\text{-red}}$ implies that $\Gamma_I (M)\in (R\text{-Mod})_{I\text{-cor}}$.  If $M\in (R\text{-Mod})_{I\text{-red}}$, then by Lemma \ref{rediscom}, $\Gamma_I(M)\in (R\text{-Mod})_{I\text{-com}}$. So, $\Gamma_I(M)\in \mathfrak{E}$.
 
 \item By Theorem \ref{TGM}, $\Lambda_I (M)\in (R\text{-Mod})_{I\text{-red}}$ for any $M\in (R\text{-Mod})_{I\text{-cor}}$. If $M\in (R\text{-Mod})_{I\text{-cor}}$, then by Lemma \ref{rediscom}, $\Lambda_I(M)\in (R\text{-Mod})_{I\text{-tor}}$. So, $\Lambda_I(M)\in \mathfrak{D}$.
 
 \item  If $M\in \mathfrak{D}$, then $\Gamma_I(M)$ is  $I$-complete (Lemma \ref{rediscom}).
 So $\Lambda_I(\Gamma_I(M))\cong \Gamma_I(M)\cong M$ since by hypothesis $M$ is $I$-torsion. On the other hand, if $M\in \mathfrak{E}$, then $\Lambda_I(M)$ is  $I$-torsion (Lemma \ref{rediscom}).
 So $\Gamma_I(\Lambda_I(M))\cong \Lambda_I(M)\cong M$ since by hypothesis $M$ is $I$-complete.
 
 \item  If $I^2=I$, then $R\text{-Mod}= (R\text{-Mod})_{I\text{-red}}= (R\text{-Mod})_{I\text{-cor}}$. So, $\mathfrak{E}=(R\text{-Mod})_{I\text{-com}}$ and $\mathfrak{D}=(R\text{-Mod})_{I\text{-tor}}$. The rest follows from part 3).
\end{enumerate}
\end{prf}

 \begin{rem}\rm
  It is tempting to think that $\mathfrak{D}=\mathfrak{E}$ is a small subcategory. However, this is not the case. For any module $M$, the $R$-modules $M/IM$ and $(0:_MI)$ belong to $\mathfrak{D}$. In particular, every module has a submodule and a quotient which is contained in $\mathfrak{D}$.   
 \end{rem}

 \subsection{Comparison with weak proregularity} 
 
 \begin{paragraph}\noi 
  The conditions: 1) weak proregularity (which is well studied in the literature) and 2) $I$-reduced and $I$-coreduced modules (being studied in this paper) are both necessary for the GM Duality and MGM Equivalence to hold. For the former, the aforementioned results hold in the derived category setting whereas for the later they hold in the module category setting. It is therefore not unreasonable to compare these two conditions. This subsection aims at achieving this.
 \end{paragraph}

 \begin{paragraph}\noi 
 Let ${\bf r} = (r_1, \cdots, r_n)$ be a sequence of elements of a ring $R$. To this sequence, we associate the Koszul complex $K(R; {\bf r})$. 
For each $i\geq 1$,  let ${\bf r}^i$ be the sequence $(r_1^i, \cdots, r_n^i)$. There is a corresponding Koszul
complex $K(R; {\bf r}^i)$. Recall that an inverse system of $R$-modules $\{M_i\}_{i\geq 1}$ is called {\it pro-zero} if for every $i$
there is some $j\geq i$ such that the $R$-homomorphism $M_j\rightarrow M_i$ is zero. 

\end{paragraph}

\begin{defn}\rm 
 A finite sequence ${\bf r} = (r_1, \cdots, r_n)$ in a ring $R$ is {\it weakly proregular} if for
every $p <0$ the inverse system of $R$-modules $\{H^p(K(R; {\bf r}^i))\}_{i\geq 1}$ is  pro-zero. 
\end{defn}
 
 \begin{defn}\rm  
  An ideal is {\it weakly proregular}  if it  is generated by a weakly proregular sequence. 
 \end{defn}
\begin{paragraph}\noi 
 If $I$ is an idempotent ideal of a ring $R$, we know that every $R$-module is both $I$-reduced and $I$-coreduced. However, $I$ need not be weakly proregular. On the other hand, every ideal $I$ of a Noetherian ring $R$ is weakly proregular but not all $R$-modules in this case are either $I$-reduced or $I$-coreduced. Instead, if the $R$-modules are finitely generated and therefore also Noetherian, then  there exists a positive integer $k$ such that all  such $R$-modules are $I^k$-reduced.
\end{paragraph}

\begin{paragraph}\noi 
Motivated by \cite{Auslander} where strongly idempotent ideals were defined for an Artin algebra, we define strongly idempotent ideals for an arbitrary ring. This definition also appears in \cite{Vyas2}.
\end{paragraph}

\begin{defn}\rm 
 An ideal $I$ of a ring $R$ is {\it strongly idempotent} if for every $i\geq 1$, $$\text{Tor}_i^R(R/I, R/I)=0.$$
\end{defn}

\begin{paragraph}\noi 
 Note that an ideal $I$ of a ring $R$ is idempotent precisely when  $\text{Tor}_1^R(R/I, R/I)=0.$ It then follows that, a strongly idempotent ideal is idempotent.
\end{paragraph}

\begin{propn}\label{ssp}
 Let $I$ be an idempotent ideal generated by a finite sequence  in a ring  $R$.  $I$  is strongly idempotent if and
only if it is weakly proregular.
\end{propn}

\begin{prf} By \cite[Proposition 4.10]{Vyas2}, an idempotent ideal $I$ is strongly idempotent if and only if the associated torsion class $\mathcal{T}_I$ is weakly stable.  However, by \cite[Theorem 4.13]{Vyas},  $\mathcal{T}_I$ is weakly stable if and only if $I$ is weakly proregular.
\end{prf}

\section{Further representability and applications}

    \begin{paragraph}\noi 
  Following \cite{Peter}, let $E_R(M)$ denote the injective hull of an $R$-module $M$ and $E$ be the direct product of the injective hulls $E_R(R/\mathfrak{m})$ of the $R$-modules $R/\mathfrak{m}$, where $\mathfrak{m}$ runs through the set of maximal ideals of $R$. $E$ is an injective cogenerator of $R$-Mod. The general Matlis duality functor is given by $$ D(-):=\text{Hom}_R(-, E).$$
  \end{paragraph}\noi

  \begin{propn}\label{L1}
   Let $R$ be a Noetherian ring and $M$ be an $R$-module.  $M$ is  reduced (resp. coreduced) if and only if the $R$-module $D(M)$ is coreduced (resp. reduced).
    \end{propn}
  
  \begin{prf}
   Suppose that $M$ is reduced. By \cite[Lemma 1.4.6]{Peter}, for any ideal $I$ of  $R$, $R/I\otimes \text{Hom}_R(M,E)\cong \text{Hom}_R(\text{Hom}_R(R/I, M), E)\cong \text{Hom}_R(\text{Hom}_R(R/I^2, M), E)\cong R/I^2\otimes \text{Hom}_R(M,E)$ so that $D(M):= \text{Hom}_R(M,E)$ is coreduced. For the converse, suppose that the $R$-module $D(M)$ is coreduced.
   By \cite[Lemma 1.4.6]{Peter},\\
   $\text{Hom}_R(\text{Hom}_R(R/I, M), E)\cong R/I\otimes \text{Hom}_R(M, E)\cong R/I^2 \otimes \text{Hom}_R(M, E)\cong$ \\  
   $\text{Hom}_R(\text{Hom}_R(R/I^2, M), E)$ for every ideal  $I$ of $R$. Since $E$ is an injective cogenerator, we have  for every ideal $I$ of $R$, $\text{Hom}_R(R/I, M)\cong \text{Hom}_R(R/I^2, M)$ and this shows that $M$  is reduced.  The second part is immediate from Proposition \ref{P1}.  
   \end{prf}

   \begin{exam}\rm Let $(R, \mathfrak{m})$ be a Gorenstein local ring of dimension $d$ and $M$ be a finitely generated $R$-module. For any $0\leq i \leq d$, the local cohomology $R$-module $H_{\mathfrak{m}}(M)$ is reduced (resp. coreduced) if and only if the $R$-module $\text{Ext}_R^{d-i}(M, R)$ is coreduced (resp. reduced). This is immediate from the local duality theorem and Proposition \ref{L1}.
    
   \end{exam}

 \subsection{Representability of $\Gamma_I$}
  
 \begin{propn}\label{P4} Let   $I$ be any  ideal of a ring $R$.
  For any two $R$-modules $M$ and $N$, where $M$ is $I$-coreduced we have 
    $$\Gamma_I(\text{Hom}_R(M,N))\cong \text{Hom}_R(\Lambda_I(M), N) ~\text{and}~  \Gamma_I(D(M)) 
    \cong D(\Lambda_I(M)).$$
  \end{propn}

 \begin{prf}  For any $I$-coreduced $R$-module $M$: 1) $\Lambda_I(M)\cong R/I\otimes M$ and;  2) the $R$-module $\text{Hom}_R(M, N)$ is $I$-reduced and therefore $\text{Hom}_R(R/I, \text{Hom}_R(M, N))\cong \Gamma_I(\text{Hom}_R(M, N))$. This together with the Hom-Tensor adjunction, we get 
 
 $$\text{Hom}_R(\Lambda_I(M), N)\cong \text{Hom}_R(R/I\otimes M, N)
 \cong  \text{Hom}_R(R/I, \text{Hom}_R(M, N))\cong $$
 $$ \Gamma_I(\text{Hom}_R(M, N)).$$
  By taking $N=E$, the second isomorphism is obtained.
 \end{prf}

  \begin{cor}\label{4.2}
   Let  $I$ be any ideal  of a ring $R$, $M$ be an $I$-coreduced $R$-module and $N$ any $R$-module. 
   
   \begin{enumerate}
    \item The $R$-module $\text{Hom}_R(\Lambda_I(M), N)$ is $I$-torsion.
    
    \item If $M$ is $I$-complete, then $\text{Hom}_R(M, N)$ (and hence $\mathcal{D}(M)$) is $I$-torsion.
    
    \item If $I^2=I$, then          
    $\Gamma_I(N)\cong \text{Hom}_R(\Lambda_I(R), N)$  and $\Gamma_I(E)=\mathcal{D}(\Lambda_I(R))$.

    \item If $I^2=I$ and $R$ is $I$-complete (resp. $\Lambda_I(R)=0$), then $N$ is $I$-torsion (resp. $I$-torsion-free).
   \end{enumerate}

  \end{cor}

  \begin{prf}
  
  \begin{enumerate}
  
\item
   By Proposition \ref{P4}, $\Gamma_I(\text{Hom}_R(M, N))\cong \text{Hom}_R(\Lambda_I(M), N)$. Taking $M=\Lambda_I(M)$   and the fact that the functor   $\Lambda_I(-)$ is idempotent on $I$-coreduced modules, we get $\Gamma_I(\text{Hom}_R(\Lambda_I(M), N))\cong \text{Hom}_R(\Lambda_I(M), N)$.     
   
    \item Immediate from  part 1.
    
    \item If $I$ is idempotent, then $R$ is an $I$-coreduced $R$-module. By Proposition \ref{P4},  $\Gamma_I(N)\cong \Gamma_I(\text{Hom}_R(R, N))\cong \text{Hom}_R(\Lambda_I(R), N)$.

    \item This is immediate from part 3.
  \end{enumerate}
  \end{prf}

  \subsection{Representability of $\Lambda_I$}
  
  \begin{propn}\label{P5}  
  Let  $I$ be a finitely generated ideal  of a ring $R$, $M$ an $I$-reduced $R$-module and  $N$  an injective $R$-module.   We have   
   $$\Lambda_I(\text{Hom}_R(M,N))\cong \text{Hom}_R(\Gamma_I(M),N)~ \text{and} ~\Lambda_I(D(M))\cong D(\Gamma_I(M)).$$
  \end{propn}

\begin{prf}
    By the fact that  $\Gamma_I(M)\cong \text{Hom}_R(R/I, M)$ and \cite[Lemma 1.4.6]{Peter}, we have  
  $\text{Hom}_R(\Gamma_I(M),N)\cong \text{Hom}_R(\text{Hom}_R(R/I, M), N)\cong R/I\otimes \text{Hom}_R(M, N) \cong \Lambda_I(\text{Hom}_R(M,N))$ since $\text{Hom}_R(M, N)$ is $I$-coreduced under the conditions given in the hypothesis.    Let $N:=E$; we get $\Lambda_I(\text{Hom}_R(M,E))\cong \text{Hom}_R(\Gamma_I(M),E)$. So that, $\Lambda_I(D(M))\cong D(\Gamma_I(M))$.
  \end{prf}

  \begin{cor}\label{4.3} 
   Let  $I$ be a  finitely generated ideal  of a ring $R$ and $N$ an injective $R$-module. If $M$ is an $I$-reduced $R$-module, then  
   
   \begin{enumerate}
    \item  the $R$-modules $\text{Hom}_R(\Gamma_I(M), N)$  and $\mathcal{D}(\Gamma_I(M))$ are $I$-complete.
\item $M$ $I$-torsion implies that $\text{Hom}_R(M, N)$ (and hence $\mathcal{D}(M)$) is $I$-complete.
      
    \item If $R$ is an $I$-reduced $R$-module,  $\Lambda_I(N)\cong \text{Hom}_R(\Gamma_I(R), N)$ and $\Lambda_I(E)\cong\mathcal{D}(\Gamma_I(R))$.
    
    \item If $R$ is $I$-reduced as an $R$-module and $I$-torsion (resp. $\Gamma_I(R)=0$), then $N$ is $I$-complete (resp. $\Lambda_I(N)=0$).
   \end{enumerate}

  \end{cor}

  \begin{prf}
  
  \begin{enumerate}
  
\item
   By Proposition   \ref{P5},  $  \Lambda_I(\text{Hom}_R(M, N))\cong \text{Hom}_R(\Gamma_I(M), N).$ Taking   $M$ to be the module $\Gamma_I(M)$  and the fact that the functor    $\Gamma_I(-)$ is idempotent, we get
   $\Lambda_I(\text{Hom}_R(\Gamma_I(M), N))\cong \text{Hom}_R(\Gamma_I(M), N)$.  
   
   \item Immediate from part 1).
   
    \item   By Proposition   \ref{P5}, $\Lambda_I(N)\cong \Lambda_I(\text{Hom}_R(R, N))\cong \text{Hom}_R(\Gamma_I(R), N)$.

   \item  Immediate from part 3).   
     
  \end{enumerate}
  \end{prf}
 
\begin{propn}
 If $I$ is an idempotent  finitely generated ideal of a ring $R$, then the functor $\Lambda_I$ is exact on a full subcategory of injective $R$-modules.
\end{propn}

\begin{prf}
 By Corollary \ref{4.3}, if $M$ is an injective $R$-module, then\\ $\Lambda_I(M)\cong \text{Hom}_R(\Gamma_I(R),M)$, i.e., on a full subcategory of injective $R$-modules the functors $\Lambda_I(-)$ and $\text{Hom}_R(\Gamma_I(R), -)$ are isomorphic. So, $\Lambda_I$ is left exact since it is isomorphic to a left exact functor $\text{Hom}_R(\Gamma_I(R), -)$. However, it is also true that since $I^2=I$, $\Lambda_I(-)\cong R/I\otimes -$ which is right exact.
\end{prf}

\begin{cor}
 For any finitely generated ideal $I$ of a ring $R$, the diagram in Figure \ref{figure2}  is commutative. In this diagram, $\mathcal{A}_I$ (resp. $\mathcal{B}_I$) denotes the full subcategory of $R$-Mod consisting of $I$-reduced (resp. $I$-coreduced) $R$-modules.
 \end{cor}
  
  \begin{figure}[!ht] 
\begin{tikzpicture}[scale=1]

\node  at (0, 0) {$\cdots$};
\draw[->, line width = 0.2mm] (0.5, 0) -- (2.5, 0);

\node  at (3, 0) {$\mathcal{A}_I$};
\draw[->, line width = 0.2mm] (3.5, 0) -- (5.5, 0);

\node  at (6, 0) {$\mathcal{B}_I$};
\draw[->, line width = 0.2mm] (6.5, 0) -- (8.5, 0);

\node  at (9, 0) {$\mathcal{A}_I$};
\draw[->, line width = 0.2mm] (9.5, 0) -- (11.5, 0);

\node  at (12, 0) {$\mathcal{B}_I$};
\draw[->, line width = 0.2mm] (12.5, 0) -- (13.5, 0);
\node  at (14, 0) {$\cdots$};
\node  at (4.4, 0.3) {$D(-)$};
\node  at (7.4, 0.3) {$D(-)$};
\node  at (10.4, 0.3) {$D(-)$};

\node  at (0, -3) {$\cdots$};
\draw[->, line width = 0.2mm] (0.5, -3) -- (2.5, -3);

\node  at (3, -3) {$\mathcal{B}_I$};
\draw[->, line width = 0.2mm] (3.5, -3) -- (5.5, -3);
\node  at (6, -3) {$\mathcal{A}_I$};
\draw[->, line width = 0.2mm] (6.5, -3) -- (8.5, -3);

\node  at (9, -3) {$\mathcal{B}_I$};
\draw[->, line width = 0.2mm] (9.5, -3) -- (11.5, -3);

\node  at (12, -3) {$\mathcal{A}_I$};
\draw[->, line width = 0.2mm] (12.5, -3) -- (13.5, -3);
\node  at (14, -3) {$\cdots$};

\node  at (4.4, -3.3) {$D(-)$};
\node  at (7.4, -3.3) {$D(-)$};
\node  at (10.4, -3.3) {$D(-)$};

\draw[->, line width = 0.2mm] (12, -0.5) -- (12, -2.5);

\draw[->, line width = 0.2mm] (9, -0.5) -- (9, -2.5);

\draw[->, line width = 0.2mm] (6, -0.5) -- (6, -2.5);

\draw[->, line width = 0.2mm] (3, -0.5) -- (3, -2.5);

\node  at (3.4, -1.5) {$\Gamma_I$};

\node  at (6.4, -1.5) {$\Lambda_I$};

\node  at (9.4, -1.5) {$\Gamma_I$};

\node  at (12.4, -1.5) {$\Lambda_I$};

\end{tikzpicture}
\caption{Commutative diagram I}\label{figure2}   
  \end{figure}

 \begin{prf}
  By Proposition \ref{P4}, $\Gamma_I(D(M))\cong D(\Lambda_I(M))$. This establishes the first commutative square. The second commutative square is established by Proposition \ref{P5} which asserts that  $\Lambda_I(D(M))\cong D(\Gamma_I(M))$. Repeating this process leads to the required commutative diagram.
 \end{prf}

\begin{cor}
 For any finitely generated idempotent ideal $I$  of a ring $R$ and for any $R$-module $M$, we have 
 \begin{enumerate}
  \item $\text{Hom}_R(R/I, \mathcal{D}(M))\cong \mathcal{D}(R/I\otimes M)$. 
  
  \item $R/I\otimes \mathcal{D}(M)\cong \mathcal{D}(\text{Hom}_R(R/I, M))$. 
  
  \item 
  The following diagram in Figure \ref{figure3} is commutative.
  \begin{figure}[!ht] 
\begin{tikzpicture}[scale=1]

\node  at (0, 0) {$\cdots$};
\draw[->, line width = 0.2mm] (0.5, 0) -- (2.3, 0);

\node  at (3, 0) {$R\text{-Mod}$};
\draw[->, line width = 0.2mm] (3.7, 0) -- (5.3, 0);

\node  at (6, 0) {$R\text{-Mod}$};
\draw[->, line width = 0.2mm] (6.7, 0) -- (8.3, 0);

\node  at (9, 0) {$R\text{-Mod}$};
\draw[->, line width = 0.2mm] (9.7, 0) -- (11.3, 0);

\node  at (12, 0) {$R\text{-Mod}$};
\draw[->, line width = 0.2mm] (12.7, 0) -- (13.3, 0);
\node  at (14, 0) {$\cdots$};

\node  at (4.4, 0.3) {$D(-)$};
\node  at (7.4, 0.3) {$D(-)$};
\node  at (10.4, 0.3) {$D(-)$};

\node  at (0, -3) {$\cdots$};
\draw[->, line width = 0.2mm] (0.7, -3) -- (2.3, -3);

\node  at (3, -3) {$R\text{-Mod}$};
\draw[->, line width = 0.2mm] (3.7, -3) -- (5.3, -3);
\node  at (6, -3) {$R\text{-Mod}$};
\draw[->, line width = 0.2mm] (6.7, -3) -- (8.3, -3);

\node  at (9, -3) {$R\text{-Mod}$};
\draw[->, line width = 0.2mm] (9.7, -3) -- (11.3, -3);

\node  at (12, -3) {$R\text{-Mod}$};
\draw[->, line width = 0.2mm] (12.7, -3) -- (13.3, -3);
\node  at (14, -3) {$\cdots$};

\node  at (4.4, -3.3) {$D(-)$};
\node  at (7.4, -3.3) {$D(-)$};
\node  at (10.4, -3.3) {$D(-)$};

\draw[->, line width = 0.2mm] (12, -0.5) -- (12, -2.5);

\draw[->, line width = 0.2mm] (9, -0.5) -- (9, -2.5);

\draw[->, line width = 0.2mm] (6, -0.5) -- (6, -2.5);

\draw[->, line width = 0.2mm] (3, -0.5) -- (3, -2.5);

\node  at (3.8, -1.5) {$\text{Hom}_R(R/I, -)$};

\node  at (6.8, -1.5) {$R/I\otimes -$};

\node  at (9.8, -1.5) {$\text{Hom}_R(R/I, -)$};

\node  at (12.8, -1.5) {$R/I\otimes -$};
\end{tikzpicture}
\caption{Commutative diagram II}\label{figure3}   
  \end{figure} 
  \end{enumerate}
\end{cor}

\subsection{Computation of natural transformations}

\begin{paragraph}\noi 
In this subsection, we demonstrate that representability of $\Gamma_I$ and $\Lambda_I$ which is facilitated by $I$-reduced and $I$-coreduced modules makes it easy (by use of Yoneda Lemma) to compute natural transformations from the functors $\Gamma_I$ and $\Lambda_I$ to other functors under some suitable conditions.
\end{paragraph}

\begin{propn}
 For any ideal $I$ of a ring $R$, and functors $$\Gamma_I: (R\text{-Mod})_{I\text{-red}} \rightarrow (R\text{-Mod})_{I\text{-cor}} ~~~\text{and}~~~ I\otimes - : (R\text{-Mod})_{I{\bf \text{-red}}}\rightarrow {\text{{\bf Set}}},$$ we  have 
 $$\text{Nat}(\Gamma_I(-), \Gamma_I(-))\cong  R/I ~~\text{and}~~\text{Nat}\left(\Gamma_I(-), I\otimes - \right)\cong 0.$$
 \end{propn}

 \begin{prf}Recall that $R/I\in(R\text{-Mod})_{I\text{-red}}$.  
  Since the functor $\Gamma_I$ is representable on $(R\text{-Mod})_{I\text{-red}}$, Yoneda Lemma asserts that for any functor $F:(R\text{-Mod})_{I\text{-red}}\rightarrow \text{{\bf Set}}$, $$\text{Nat}(\text{Hom}_R(R/I, -), F(-))\cong F(R/I).$$ It follows that $\text{Nat}(\Gamma_I(-), \Gamma_I(-))\cong \Gamma_I(R/I)\cong (0_{R/I}I)\cong R/I.$
 \end{prf}

 \begin{paragraph}\noi 
  Let $R$ be a finite dimensional algebra over a field $k$ and $R$-mod a category of all finitely generated $R$-modules. Let $\mathcal{D}:=\text{Hom}_k(-, k)$. The Nakayama functor $\mathcal{V}$ is the composition of two contravariant functors $\text{Hom}_R(-, R)$ and $\mathcal{D}$, i.e.,  $$\mathcal{V}:= \mathcal{D}\text{Hom}_R(-, R): R\text{-mod}\rightarrow R\text{-mod}.$$ 
 If $R$-proj (resp. $R$-inj) denotes the full subcategory of $R$-mod consisting of all projective $R$-modules (resp. injective $R$-modules), then the restriction of $\mathcal{V}$ on $R$-proj defines an equivalence  
 $$\mathcal{V} : R\text{-proj}\rightarrow R\text{-inj}$$ between $R\text{-proj}$ and $R\text{-inj}$ with quasi-inverse $\mathcal{V}^{-1}:=\text{Hom}_R(\mathcal{D}(R), -) $, see \cite[Definition 2.8 \& Proposition 2.10]{Assem}.
  \end{paragraph}
 
 \begin{propn}Let $R$ be a Noetherian von-Neumann regular finite dimensional algebra over a field $k$. If $R$ is self-injective, then for any $M\in R\text{-inj}$ and  functors $$\Lambda_I : R\text{-inj}\rightarrow \text{\bf Set}~~ \text{and}~~\mathcal{V}^{-1} : R\text{-inj}\rightarrow R \text{-proj}, $$ we have $$\text{Nat}\left(\Lambda_I(M), \mathcal{V}^{-1}(M)\right)\cong k.$$
 \end{propn}
 
 \begin{prf}
  By hypothesis, $R\in R\text{-inj}$. By \cite{Brodmann}$,\Gamma_I$ preserves injective modules defined over Noetherian rings. So, $\Gamma_I(R)\in R\text{-inj}$. By Corollary \ref{4.3}(3), $\Lambda_I(M)\cong \text{Hom}_R(\Gamma_I(R), M)$. Since $R$ is Noetherian, there exists a positive integer $t$ such that $\Gamma_I(R)\cong \text{Hom}_R(R/I^t, R)$. We can see that the functor $\Lambda_I : R\text{-inj}\rightarrow \text{\bf Set}$ is representable. So, given the functor $\mathcal{V}^{-1} : R\text{-inj}\rightarrow R \text{-proj}$, we invoke Yoneda Lemma to have 
$\text{Nat}\left(\Lambda_I(M), \mathcal{V}^{-1}(M)\right)\cong \text{Nat}\left(\text{Hom}_R(\Gamma_I(R), M), \mathcal{V}^{-1}(M)\right)\cong  \mathcal{V}^{-1}(\Gamma_I(R))$. However, $\mathcal{V}^{-1}(\Gamma_I(R))\cong \text{Hom}_R(\mathcal{D}(R), \Gamma_I(R))$. Substituting for $\mathcal{D}(R)$ and $\Gamma_I(R)$; and applying Yoneda embedding,  we get 
$$\mathcal{V}^{-1}(\Gamma_I(R))\cong \text{Hom}_R(\text{Hom}_k(R, k), \text{Hom}_R(R/I^t, R))\cong \text{Hom}_k(k, R)\cong k.$$
  \end{prf}
  
  \begin{lem} \label{lm} Let $I$ be an ideal of a ring $R$. If $R$ is a reduced (resp. $I$-reduced) $R$-module, then the full subcategory $R$-proj of  $R$-Mod which consists of all projective $R$-modules has all  modules reduced (resp. $I$-reduced).  
  \end{lem}

\begin{prf}
 If $R$ is $I$-reduced, then so is any free $R$-module; since such a module is isomorphic to $R^n$ for some positive integer $n$ and $I$-reduced modules are closed under taking direct sums. Since every projective $R$-module is a direct summand of a free $R$-module, it follows that, all projective $R$-modules in this case are $I$-reduced. For the reduced case, the proof is similar.
\end{prf}

\begin{propn}
 Let $I$ be an ideal of a ring $R$ such that $R$ is an  $I$-reduced $R$-module and  $R/I$ is a projective $R$-module. For any  $M\in R\text{-proj}$ and functors: $$ \Gamma_I : R\text{-proj}\rightarrow {\text{\bf Set}}~~\text{and} ~~ \mathcal{V} : R\text{-proj}\rightarrow R\text{-inj};$$
 
 $$\text{Nat}(\Gamma_I(M), \mathcal{V}(M))\cong \mathcal{D}(0:_RI).$$
\end{propn}

 \begin{prf} If $M\in R\text{-proj}$, then by Lemma \ref{lm}, it is $I$-reduced and therefore $\Gamma_I(M)\cong \text{Hom}_R(R/I, M)$. The hypothesis satisfies conditions of the Yoneda Lemma. So, we have 
 $$\text{Nat}(\Gamma_I(M), \mathcal{V}(M))\cong  \text{Nat}(\text{Hom}_R(R/I,  M), \mathcal{V}(M))\cong \mathcal{V}(R/I).$$ However, $$ \mathcal{V}(R/I)\cong \mathcal{D}\text{Hom}_R(R/I, R)\cong \mathcal{D}(0:_RI).$$
  
 \end{prf}
 Data sharing not applicable to this article as no datasets were generated or analysed during the current study.

\section*{Funding and/or Competing interests.}
 
 \begin{paragraph}\noi
 The author was supported by the International Science Programme through the Eastern Africa Algebra Research Group and also by the EPSRC GCRF project EP/T001968/1, Capacity building in Africa via technology-driven research in algebraic and arithmetic geometry (part of the Abram Gannibal Project). Part of this work was written while the author was visiting University of Glasgow and University of Oxford. He is grateful to Kobi Kremnizer, Balazs Szendroi   and Michael Wemyss for the stimulating discussions  and for the hospitality during his stay in the United Kingdom. To Annet Kyomuhangi and Amanuel Senbato Mamo, thank you for all the comments. 
 \end{paragraph}


\addcontentsline{toc}{chapter}{Bibliography}
\begin{thebibliography}{99}


\bibitem{Alonso}
L. Alonso Tarrı́o, A. Jeremı́as López and J. Lipman, Local homology and cohomology on schemes, {\it Annales Scientifiques de l'École Normale Supérieure}, (4), {\bf 30}(1),  (1997), 1--39.

\bibitem{Assem} I. Assem, D. Simson, and A. Skowronski,  ``Elements of the Representation Theory of Associative Algebras: Volume 1: Techniques of Representation Theory'', Cambridge University Press, 2006.

\bibitem{Auslander}  M. Auslander, M. I. Platzeck and  G. Todorov, Homological theory of idempotent ideals. {\it Trans. Amer. Math. Soc.} {\bf 332}(2),  (1992), 667–692.
 

\bibitem{Barthel} T.  Barthel, D.  Heard and G. Valenzuela, Local duality in algebra and topology, {\it Adv. Math.}  {\bf 335},  (2018), 563-663.

\bibitem{Barthel2} T.  Barthel, D.  Heard and G.  Valenzuela, Derived completion for comodules, {\it Manuscripta Math.}. {\bf 161}, (2020), 409-438. 

\bibitem{Brodmann} 
M. Brodmann and R. Y. Sharp, ``Local cohomology: An algebraic introduction with geometric applications'',
Cambridge Studies in Advanced Mathematics, 136, Cambridge University Press, Cambridge, Second Edition, 2013.
     
 \bibitem{Divaani}    
   K. Divaani-Aazar,  Local homology, finiteness of Tor modules and cofiniteness, {\it J. Algebra Appl.}
16(12), (2017),  1750240 (10 pages)
 World Scientific Publishing Company
DOI: 10.1142/S0219498817502401.

\bibitem{glduality} T. H. Freitas, V. H. Jorge-Pérez, C. B. Miranda-Neto and  P. Schenzel, Generalized local duality, canonical modules, and prescribed bound on projective dimension, https://arxiv.org/abs/2112.12632.

 \bibitem{Greenless} J. P. C. Greenless and J. P. May, Derived functors of $I$-adic completion and local  homology, {\it J. Algebra}, {\bf 149}, (1992), 438--453.
 
 
 \bibitem{Hartshorne}R. Hartshorne,   ``Local cohomology'', Lecture notes in mathematics, 1967. A seminar given by A. Grothendieck, Harvard University. Fall,  (Vol. 41). Springer Verlag, Berlin-New York, 1961.
 
 
\bibitem{Kempf} G. Kempf, The Grothendieck-Cousin complex of an induced representation, {\it Adv.  Math.}, {\bf 29},  (1978), 310--396.
 
\bibitem{Kyom} A. Kyomuhangi and D. Ssevviiri, The locally nilradical for modules over commutative rings, {\it Beitr. Algebra Geom.}, {\bf 61}, (2020), 759–769.


\bibitem{Kyom3} A. Kyomuhangi and D. Ssevviiri, Generalised reduced modules, {\it Rend. Circ. Mat. Palermo (2)}, https://doi.org/10.1007/s12215-021-00686-8, (2021).
  
\bibitem{Kyom2} A. Kyomuhangi and D. Ssevviiri, Coreduced modules and the associated radicals, submitted.
 
\bibitem{Lee}
T. K. Lee,   and Y. Zhou, ``Reduced Modules'', Rings, Modules, Algebras and Abelian Groups. Lecture Notes in
Pure and Applied Math, vol. 236, pp. 365--377. Marcel Dekker, New York, 2004.

\bibitem{Mac} S. MacLane, ``Categories for the working mathematician'', Berlin, Heidelberg, New York: Springer-Verlag, 1971.


\bibitem{Yekutieli}  M.  Porta,  L. Shaul and A. Yekutieli, On the homology of completion and torsion, {\it Algebr. Represent. Theory}, {\bf 17}, (2014),  31–67.

\bibitem{Rege} M. B. Rege and A. M.  Buhphang,  On reduced modules and rings. {\it Int. Electron. J. Algebra} {\bf 3}, (2008), 58–74.

\bibitem{Peter}  P. Schenzel and A. Simon, ``Completion, Cech and Local Homology and Cohomology'', Springer Monographs in Mathematics, 2018.
 

\bibitem{Sharp} R. Y. Sharp. Local cohomology theory in commutative algebra. {\it Q. J. Math.}, {\bf 21} (4),   (1970), 425–434.

\bibitem{David}  D. Ssevviiri, Nilpotent elements control the structure of a module. arXiv:1812.04320 [math.RA],  (2018).

\bibitem{Suzuki}N. Suzuki, On the generalized local cohomology and its duality, {\it J. Math. Kyoto Univ}. {\bf 18}, (1978),  71–85.

\bibitem{Vyas} R. Vyas and  A. Yekutieli, Weak proregularity, weak stability, and the noncommutative MGM equivalence. {\it J. Algebra},
{\bf 513}, (2018), 265–325.

\bibitem{Vyas2} R. Vyas,  Weakly Stable torsion classes, {\it Algebr. Represent. Theory}, {\bf 22}, (2019),  1183--1207.  

\bibitem{Weibel} C. A. Weibel, ``An introduction to homological algebra'', Cambridge studies in advanced mathematics 38, Cambridge University Press, 1994.

\bibitem{Yassemi} S. Yassemi, Generalized section functors, {\it J. Pure Appl. Algebra}, {\bf 95},  (1994), 103–119.

 
\bibitem{Y2} A. Yekutieli,  On flatness and completion for infinitely generated modules over Noetherian rings, {\it Comm. Algebra}, {\bf 39}(11),  (2011), 4221--4245.

\end{thebibliography}
\end{document}